\documentclass[12pt]{article}
\usepackage{amsmath, amssymb,amsgen, amsthm, amscd}
\usepackage[all]{xy}
\usepackage{hyperref}

\newtheorem{thm}{Theorem}
\newtheorem{lemma}[thm]{Lemma}
\newtheorem{cor}[thm]{Corollary}
\newtheorem{prop}[thm]{Proposition}

\theoremstyle{definition}

\newtheorem{rem}[thm]{Remark}

\theoremstyle{remark}

\newcommand{\ca}{\mathcal{A}}
\newcommand{\cb}{\mathcal{B}}
\newcommand{\cc}{\mathcal{C}}
\newcommand{\cl}{\mathcal{L}}
\newcommand{\cg}{\mathcal{G}}

\newcommand{\dint}{{\displaystyle \int \!\!\!\!\!\!-}}

\newcommand{\ch}{\mathcal{H}}

\newcommand{\al}{\alpha}
\newcommand{\lb}{\lambda}

\newcommand{\de}{\delta}
\newcommand{\De}{\Delta}

\newcommand{\w}{\omega}
\newcommand{\e}{\varepsilon}

\newcommand{\W}{\Omega}

\newcommand{\g}{\gamma}
\newcommand{\ve}{\varepsilon}
\newcommand{\G}{\Gamma}

\newcommand{\s}{\sigma}

\newcommand{\te}{\theta}

\newcommand{\del}{\partial}

\newcommand{\ify}{\infty}

\newcommand{\bg}{\mathrm{B}G}
\newcommand{\eg}{\rm{E}G}

\newcommand{\Cb}{\mathbb{C}}
\newcommand{\Rb}{\mathbb{R}}

\newcommand{\Zb}{\mathbb{Z}}

\DeclareMathOperator{\Diff}{Diff}

\DeclareMathOperator{\tr}{tr}

\DeclareMathOperator{\Dom}{Dom}
\DeclareMathOperator{\Ran}{Ran}
\DeclareMathOperator{\supp}{supp}

\newcommand{\chg}{\ch\left(GL_n(\Rb)\right)}
\newcommand{\grch}{\ch^*}

\newcommand{\rp}{\rho^+}
\newcommand{\rn}{\rho^-}
\newcommand{\hz}{h_{(0)}}
\newcommand{\ho}{h_{(1)}}
\newcommand{\htw}{h_{(2)}}
\newcommand{\hth}{h_{(3)}}
\newcommand{\hf}{h_{(4)}}
\newcommand{\hfi}{h_{(5)}}

\newcommand{\gz}{g_{(0)}}
\newcommand{\go}{g_{(1)}}
\newcommand{\gtw}{g_{(2)}}
\newcommand{\gth}{g_{(3)}}

\begin{document}

\title{
Secondary Characteristic Classes and
Cyclic Cohomology of Hopf Algebras
}
\author{Alexander Gorokhovsky
\\Department of Mathematics,\\
 University of Michigan,\\
Ann Arbor, MI 48109-1109, USA\\
gorokhov@math.lsa.umich.edu}
\maketitle

\begin{abstract}

Let $X$ be a manifold on which discrete (pseudo)group
  of diffeomorphisms $\Gamma$
acts, and let $E$ be a $\Gamma$-equivariant vector bundle on $X$.
We give a construction of cyclic cocycles on the cross product algebra
$C_0^{\infty}(X) \rtimes \Gamma$ representing the
 equivariant characteristic classes of
$E$. Our formulas generalize Connes'  Godbillon-Vey cyclic cocycle.
An essential tool of our construction is 
 Connes-Moscovici's theory of cyclic
cohomology of Hopf algebras.

\end{abstract}

\numberwithin{equation}{section}
\section{Introduction}
In the paper \cite{cn86} A. Connes provided an explicit construction of
the Godbillon-Vey cocycle in the cyclic cohomology. 
The goal of this
paper is to give a similar construction for the higher secondary classes.

First, let us recall Connes' construction. Let $M$ be a smooth oriented
manifold and let $\G \subset \Diff^+(M)$ be a discrete group of
orientation-preserving diffeomorphisms of $M$. Let $\w$ be a volume
form on $M$. Define the following function on $M\times \G$:
 $\de(g)=\frac{\w^g}{\w}$, where the superscript denotes the group
action. Then one can define one-parametric group of diffeomorphisms
of the algebra $\ca=C_0^{\ify}(M)\rtimes \G$  by
\begin{equation}\label{group}
\phi_t (aU_g)=a\de(g)^tU_g
\end{equation}
This is Tomita-Takesaki group of automorphisms, associated to the weight
on $\ca$ given by $\w$.

Consider now the transverse fundamental  class -- cyclic $q$-cocycle on
$\ca$ given by
\begin{multline}
\tau( a_0U_{g_0}, a_1U_{g_1},\dots, a_qU_{g_q})=\\
\begin{cases}
\frac{1}{q!}\int \limits_M a_0 da_1^{g_0}da_2^{g_0g_1}\dots
da_q^{g_0g_1\dots g_{q-1}} &\text{ if } g_0g_1\dots g_q=1\\
0 &\text{ otherwise}
\end{cases}
\end{multline}

To study the behavior of this cocycle under the 1-parametric group
\eqref{group}, consider the ``Lie derivative'' $\cl$ acting on the
cyclic complex by
$\cl \xi =\frac{d}{dt}|_{t=0}\phi_t^* \xi$, $xi$ being a cyclic cochain.
It turns out that in general $\tau$ is not invariant under the group
\eqref{group}, and $\cl \tau \ne 0$.

However, it was noted by Connes that one always has 
\begin{equation}\label{inv}
\cl^{q+1} \tau=0
\end{equation}
and $\cl^q \tau$ is invariant under the action of the group \eqref{group}.
One deduces from this that if $\iota_{\de}$ is the analogue of the
interior derivative (see \cite{cn86}), then $\iota_{\de}\cl^q\tau$ is a
cyclic cocycle.

This is Connes' Godbillon-Vey cocycle. It can be related to the
Godbillon-Vey class as follows. Let $[GV]\in H^*(M_{\G})$
be the Godbillon-Vey class, where $M_{\G}=M\times_{\G}\rm{E}\G$ is the
homotopy quotient. Connes defines canonical map $\Phi:
H^*(M_{\G})\to HP^*(\ca)$. Then one has
\begin{equation}\label{phigv}
\Phi([GV])=[\iota_{\de}\cl^q\tau]
\end{equation}

The class of this cocycle is independent of the choice of the volume
form. To prove this one  can use Connes'
noncommutative Radon-Nicodym theorem 
to conclude that if one changes the volume form, the
one-parametric group $\phi_t$ remains the same modulo inner
automorphisms. 

A natural problem then is to extend this construction to the cocycles
corresponding to the other secondary characteristic classes.
It was noted by Connes \cite{cn94} that if instead of 1-dimensional
bundle of $q$-forms on $M$ one considers $\G$ equivariant trivial
bundle of rank $n$, then in place of 1-parametric group \eqref{group} one encounters
coaction of the group $GL_n(\Rb)$ on the algebra $\ca$. The difficulty
is that for $n>1$ this group is not commutative, and one can not replace
this coaction by the action of the dual group, similarly to \eqref{group}.

In this paper we show that Connes-Moscovici theory of cyclic
cohomology
for Hopf algebras (cf \cite{cm98, cm99}) provides a natural framework for the
higher-dimensional situation and allows one to give construction of the
secondary characteristic cocycles.

The situation we consider is the following. We have an
orientation-preserving 
 action of
discrete group (or pseudogroup) $\G$ on the oriented manifold $M$, and a
trivial bundle $E$ on $M$ equivariant with respect to this action.
Well-known examples in which such a situation arises are the following
(cf. \cite{cn94}, \cite{ni96},\cite{cm98}).
Let $V$ be a manifold on which a discrete group (or pseudogroup)
$G$ acts, and let be $E_0$ a bundle (not necessarily trivial) on
$V$, equivariant with respect to the action of $G$. Let $U_i$, $i \in I$,
be an open
covering of $V$ such that restriction of $F$ on each $U_i$ is trivial.
Put $M=\sqcup U_i$, and let $E$ be the pull-back of $E_0$ to $M$ by the
natural projection. Then we have an action of the following pseudogroup
$\G$ on $M$: $\G=\{g_{i,j}|\ g\in G\, i,j \in I\}\cup id$,
where $\Dom g_{i,j}=g^{-1}\left(U_j\right)\cap U_i \subset U_i$,
$\Ran g_{i,j}=g\left(U_i\right)\cap U_j \subset U_j$, and the natural
composition rules. The bundle $E$ is clearly equivariant with respect to
this action. Our construction, described below, provides classes in
the cyclic cohomology of the cross-product algebra
$C_0^{\ify}(M)\rtimes \G$, rather than in the cyclic cohomology of
$C_0^{\ify}(V)\rtimes G$.  
 However, cross-product algebras $C_0^{\ify}(M)\rtimes \G$
and $C_0^{\ify}(V)
\rtimes G$ are Morita equivalent, and hence have the same cyclic
cohomology. 

Another natural example is provided by the manifold $V$ with foliation
$F$, and a bundle $E_0$ which is holonomy equivariant. We can always
choose (possibly disconnected) complete transversal $M$, such that
restriction $E$ of $E_0$ to $M$ is trivial. Let $\G$ be the holonomy
pseudogroup acting on $M$. $E$ is clearly equivariant with respect to
this action. In this case again the cross-product algebra $C_0^{\ify}(M)\rtimes
\G$ is Morita equivalent to the full algebra of the foliation
$C_0^{\ify}(V,F)$. 

We construct then a map from the cohomology of the truncated
Weil algebra (cf. e.g. \cite{kt75}) $W(\mathfrak{g}, O_n)_q$ to the periodic
cyclic cohomology $HP^*(\ca)$ of the algebra $\ca=C_0^{\ify}(M)\rtimes
\G$. The construction is the following. We consider the action of the
\emph{differential graded} Hopf algebra $\chg$ of  differential forms on
the group $GL_n(\Rb)$ on the differential graded algebra $\W^*_0(M)
\rtimes \G$, where $\W^*_0(M)$ denotes the algebra of compactly
supported differential forms on $M$. The use of differential graded
algebras allows one to conveniently encode different identities, similar
to \eqref{inv}.
We then show that Connes-Moscovici theory (or rather differential graded
version of it) allows one to define a map from the cyclic complex of
$\chg$ to the cyclic complex of $\W^*_0(M)
\rtimes \G$. We then relate cyclic complex of $\chg$ to the Weil
algebra, and cyclic cohomology of $\W^*_0(M)
\rtimes \G$ to the cyclic cohomology of $\ca$.

The paper is organized as follows.
In the next two sections we discuss cyclic complexes for differential
graded algebras and differential graded Hopf algebras respectively.
In the section \ref{S:prop} we show that two different Hopf actions,
which coincide ``modulo inner automorphisms'' induce the same
Connes-Moscovici characteristic map in cyclic cohomology, and discuss
some other properties of the characteristic map.
In the section \ref{S:forms} we construct the action of
$\chg$ on $\W^*_0(M)
\rtimes \G$. In the section \ref{S:rel} we relate cyclic complex of the
Hopf algebra $\chg$ with the Weil algebras. Finally, in the section
\ref{S:other} we prove an analogue of the formula \eqref{phigv} for the
cocycles we construct.

I would like  to thank D. Burghelea and H. Moscovici
for many helpful discussions.

\section{Cyclic complex for differential graded algebras}
In this section we collect some preliminary standard
facts about cyclic cohomology
of the differential graded algebras, and give cohomological version of
some results of \cite{goo85}.

Recall that the cyclic module $X^*$ is given by the cosimplicial module with
the face maps  $\de_i:X^{n-1}\to X^{n}$ and
and degeneracy maps $\s_i: X^n\to X^{n-1} $ $0\le i \le n$, satisfying
the usual axioms. In addition, we have for each $n$ an action of $\Zb_{n+1}$ on
$X^n$, with the generator $\tau_n$ satisfying
\begin{equation}
\tau_n  \delta_i = \delta_{i-1}  \tau_{n-1} \text{ for }1 \leq i \leq n
\text{ and }\tau_n  \delta_0 = \delta_n  
\end{equation}
\begin{equation}
\tau_n  \s_i = \s_{i-1}  \tau_{n+1} \text{ for }1 \leq i \leq n \text{
and } \tau_n  \s_0 = \s_n  \tau_{n+1}^2
\end{equation}
\begin{equation}
\tau_n^{n+1} = id
\end{equation}
For every cyclic object one can construct operators
$b: X^n \to X^{n+1}$ and $B:
X^n \to X^{n-1}$, defined by the formulas
\begin{equation}
b=\sum \limits_{j=0}^n (-1)^j d_j
\end{equation}
\begin{equation}
B=\left(\sum \limits_{j=0}^{n-1}(-1)^{j(n-1)} \tau_{n-1}^j\right) \s_{n+1} (1-(-1)^{n-1}\tau_n)
\end{equation} 
where
\begin{equation}
\s_{n+1}=\s_n\tau_{n+1}
\end{equation}
These operators satisfy
\begin{align}
b^2&=0\\
B^2&=0\\
bB+Bb&=0
\end{align}

Hence for any cyclic object $X^*$ we can construct a bicomplex
$\cb^{*,*}(X)$ as follows: $\cb^{p,q}$, $p$, $q\ge 0$ is $X^{p-q}$, or $0$ if
$p<q$, and the differential $\cb^{p,q} \to \cb^{p+1,q}$ (resp.
$\cb^{p,q+1}$) is given by $b$ (resp. $B$). By removing restriction $p$,
$q\ge 0$ we obtain periodic bicomplex $\cb_{\text{per}}$. Notice that it
has periodicity induced by the tautological shift $S:
\cb_{\text{per}}^{p,q}\to \cb_{\text{per}}^{p+1,q+1}$.

Let now $\W^*$ be a unital  graded (DG) algebra, positively graded.
We can associate with it a cyclic object
as follows (the differential $d$ is not used in this definition).

Let $C^k(\W^*)$ be the space of continuous $k+1$-linear functionals
on $\W^*$. The face and degeneracy maps are given by
\begin{align}
(\de_j\phi) (a_0,a_1,\dots,a_{k+1})=&
\phi(a_0,\dots, a_ja_{j+1},\dots,a_{k+1}) \text{ for } 0\le j \le
n-1 \nonumber \\
(\de_n \phi) (a_0,a_1,\dots,a_{k+1})=&
(-1)^{\deg a_{k+1}(\deg a_0+\dots +\deg a_k)} \phi (a_{k+1}a_0, a_1, \dots a_k)
\end{align}
\begin{equation}
(\s_j \phi)(a_0,\dots,a_{k-1})= \phi (a_0,\dots, a_{j}  , 1, a_{j+1} ,\dots a_{k-1})
\end{equation}
and the cyclic action is given by
\begin{equation}
(\tau_k \phi)(a_0,\dots,a_{k})=(-1)^{\deg a_k( \deg a_0+\dots \deg a_{k-1})}
\phi(a_k,a_0,\dots,a_{k-1})
\end{equation}
Cohomology of
total complex of the bicomplex
 $\cb$ (resp. $\cb_{\text{per}}$) where we consider only \emph{finite}
cochains, is the cyclic (resp. periodic cyclic) cohomology of $\W^*$, for
which we use notation $HC^*(\W^*)$ (resp. $HP^*(\W^*)$).

Suppose now that $\W^*$ is a differential graded (DG) algebra with the
differential of degree $1$. 
We say that $\phi \in C^k(\W^*)$ has weight $m$
if $\phi (a_0,a_1,\dots,a_k)=0$ unless $\deg a_0+\deg a_1+\dots+\deg
a_k=m$. We denote by $C^{k,p}(\W^*)\subset C^k(\W^*)$ set of weight $-p$
functionals.
Notice that in this case each $C^k(\W^*)$ is a complex in its own right,
with the grading defined above and the differential $(-1)^kd$, where 
we extend $d$ to $C^k(\W^*)$ by
\begin{equation}
d \phi(a_0,a_1,\dots,a_k)=
\sum \limits_{j=0}^k (-1)^{\deg a_0+\dots \deg a_{j-1}}
\phi( a_0,\dots, da_j,\dots, a_k)
\end{equation}

Then $db-bd=0$, $dB-Bd=0$, and
hence in this situation $\cb$ and $\cb_{\text{per}}$ become actually
tricomplexes.
Cyclic (resp. periodic cyclic) cohomology of the DG algebra $(\W^*,d)$
is then defined as the cohomology of the total complex of the tricomplex
 $\cb$ (resp. $\cb_{\text{per}}$) where we consider only \emph{finite} cochains.  
Notations for the cyclic and periodic cyclic cohomologies are
$HC^*\left((\W^*,d)\right)$ and $HP^*\left((\W^*,d)\right)$. 

One can show that cyclic cohomology can be computed by the
\emph{normalized} complex, i.e. one where cochains satisfy
\begin{equation}
\phi(a_0,a_1,\dots,a_k)=0 \text{ if } a_i=1,\ i\ge 1
\end{equation}
We will need the following result about the cyclic cohomology.
\begin{thm}
Let $\ca=\W^0$ be the $0$-degree part of $\W^*$ (which we consider as a
trivially graded algebra with the zero differential). We then have a
natural map of (total) complexes
$I:\cb_{\text{per}}(\ca) \to \cb_{\text{per}}\left((\W^*,d)\right)$
(extension of polilinear forms by $0$ to $\W^*$).
 Then the induced map in cohomology is an isomorphism.
$HP^*(\ca) \to HP^*\left((\W^*,d)\right)$. 
\end{thm}\label{isom}
To prove the theorem and to write an explicit formula for the map 
$R: \cb_{\text{per}}\left((\W^*,d)\right) \to \cb_{\text{per}}(\ca)$,
inducing the inverse isomorphism in the periodic cyclic cohomology, we
need 
the following fact (Rinehart formula) which we use
as stated in \cite{goo85}.

Let $D$ be a derivation of the graded algebra $\W^*$ of degree $\deg D$,
i.e. a linear map $D: \W^* \to \W^{*+\deg D}$ satisfying
\begin{equation}
D(ab)=(Da) b+ (-1)^{\deg D \deg a}a D(b)
\end{equation}  
It defines an operator on the complex $\cb(\W^*)$, by
\begin{equation}
\cl_D\phi(a_0, a_1,\dots,a_k)=
\sum \limits_{i=0}^k (-1)^{\deg D(a_0+\dots+a_{i-1})}
\phi(a_0,\dots,D(a_i),\dots, a_k)
\end{equation}
which commutes with $b$, $B$. The action of this operator on the periodic
cyclic
bicomplex is homotopic to zero, with the homotopy constructed as follows.
Define operators $e_D:C^{k-1}(\W^*)\to e_D:C^{k}(\W^*) $, $E_D:e_D:C^{k+1}(\W^*)
\to e_D:C^{k}(\W^*)$ by
\begin{align}
e_D\phi(a_0,a_1,\dots a_k)&= (-1)^{k+1}\phi(D(a_k)a_0,a_1\dots,a_{k-1}) \\
E_D\phi(a_0,a_1,\dots a_k)&= \sum \limits_{1\le i\le j \le k}
(-1)^{ik+1} \phi(1,a_i,a_{i+1},\dots,a_{j-1}, Da_j,\dots, a_k, a_0,\dots)
\end{align}
Then
\begin{equation}
[b+B,e_D+E_D]=\cl_D
\end{equation}

We now proceed with the proof of the Theorem \ref{isom}
\begin{proof}[Proof of the Theorem \ref{isom}]
Consider the derivation $D$ of degree $0$ given by $Da=(\deg a) a$.
On the polilinear form of weight $m$ $D$ acts by $m$. Define the
homotopy $h$ to be $\frac{1}{m}(e_D+E_D)$ on the forms of weight
$m>0$ and $0$ on  the forms of weight $0$.
We define map of complexes $R: \cb_{\text{per}}\left((\W^*,d)\right) \to
\cb_{\text{per}}(\ca)$ by
\begin{equation}\label{defR}
R \phi =c_{k,m}(dh)^m \phi \text{ for } \phi \in
C^{k,m} 
\end{equation}
where
\begin{equation}
c_{k,m}=(-1)^{km +\frac{m^2-m}{2}}
\end{equation}
This is a map of (total) complexes. Indeed, using identities
$(b+B)h+h(b+B)=id$ and $(b+B)d-d(b+B)=0$ we have:
\begin{multline}
(b+B)R\phi -R(b+B +(-1)^kd)\phi=\\
c_{k,m}
\left(( (b+B)(dh)^m-(-1)^m(dh)^m(b+B))+(-1)^{m}(dh)^{m-1}d\right)\phi=\\
c_{k,m}
\left(-(-1)^m(dh)^{m-1}d+(-1)^m(dh)^{m-1}d\right)\phi=0
 \end{multline}
 It is clear that
 \begin{equation}
 R\circ I=id
 \end{equation}
 As for $I\circ R$ we have
 \begin{equation}
 I\circ R=id -\left(\del \circ H+H\circ \del\right)
 \end{equation}
 where $\del=b+B\pm d$ -- total differential in the complex
 $\cb_{\text{per}}\left((\W^*,d)\right)$, $\pm d$ being $(-1)^kd$ on
$C^{k,m}$, and the homotopy $H$ is given
by the formula
\begin{equation}
H \phi=
\sum \limits_{j=0}^{m-1}c_{k,j}
h(dh)^j\phi \text{  for }\phi \in C^{k,m}
\end{equation}
This equality is also verified by direct computation.

Indeed, we have, for actions  on $C^{k,m}$
\begin{align}
H\circ(b+B)=\sum \limits_{j=0}^{m-1}c_{k,j}((dh)^j-(b+B)h(dh)^j)+
\sum \limits_{j=1}^{m-1}c_{k,j}(-1)^j(hd)^j \\
(b+B)\circ H=\sum \limits_{j=0}^{m-1}c_{k,j}(b+B)h(dh)^j \\
H\circ (\pm d)=\sum \limits_{j=0}^{m-2}(-1)^{k}c_{k,j}(hd)^{j+1}=
\sum \limits_{j=1}^{m-1}(-1)^{k}c_{k,j-1}(hd)^{j }\\
(\pm d)\circ H=\sum \limits_{j=0}^{m-1}(-1)^{k-j-1}c_{k,j}(dh)^{j+1}=
\sum \limits_{j=1}^{m-1}(-1)^{k-j}c_{k,j-1}(dh)^{j }
\end{align}
and adding these equalities we get the desired result.
\end{proof}

\section{Cyclic complex for differential graded Hopf algebras}
In this section we reproduce Connes-Moscovici's construction of the
cyclic module of a Hopf algebra (cf. \cite{cm98,cm99})
 in the differential graded context.

Let us start with the graded Hopf algebra $\ch^*$. We need to fix a
modular pair, i.e. a homomorphism $\de :\ch^* \to \Cb$ and a group-like
element $\s \in \ch^0$. Using the standard notations for the coproduct
and antipode, define the twisted antipode $\widetilde{S}_{\de}$ by
\begin{equation}
\widetilde{S}_{\de}(h)= \sum S(\hz)\de(\ho)
\end{equation} 
Suppose that the following condition holds:
\begin{equation}
\left(\s^{-1} \widetilde{S}_{\de}\right)^2=id
\end{equation}

Then Connes and Moscovici show that one can define a cyclic object
$\left(\grch\right)^{\sharp}= \left\{ \left(\grch \right)^{\otimes n} \right\}_{n\ge 1} $
 as follows. Face and degeneracy operators are given by
\begin{align}
\de_0 (h^1 \otimes \ldots \otimes h^{n-1})  &=
1 \otimes h^1 \otimes \ldots \otimes h^{n-1} \nonumber \\
\de_j (h^1 \otimes \ldots \otimes h^{n-1})  &= h^1 \otimes \ldots \otimes \De h^j \otimes 
\ldots \otimes h^n \text{ for }  1 \le j \le n-1, \nonumber\\
\de_n (h^1 \otimes \ldots \otimes h^{n-1})  &= h^1 \otimes \ldots \otimes h^{n-1}
\otimes \sigma \nonumber\\
\s_i (h^1 \otimes \ldots \otimes h^{n+1})  &= h^1 \otimes \ldots \otimes \ve (h^{i+1}) 
\otimes \ldots \otimes h^{n+1}  
\end{align} 
The cyclic operators are given by
\begin{multline}
\tau_n(h^1 \otimes \ldots \otimes h^{n+1})=\\
\sum (-1)^{\sum \limits_{j>i\ge 0}\deg h^1_{i} \deg h^j} \left(\widetilde{S} h^1 \right)_{(0)}h^2 
\otimes \dots \otimes \left(\widetilde{S} h^1 \right)_{(n-2)} h^n \otimes
\left(\widetilde{S} h^1 \right)_{(n-1)} \sigma
\end{multline}
where
\[
\left(\De^{n-1}\widetilde{S}h^1\right)=\sum  \left(\widetilde{S} h^1 \right)_{(0)}
\otimes \dots \left(\widetilde{S} h^1 \right)_{(n-1)}
\]
It is verified in \cite{cm99} that the operations above indeed define a
structure of a cyclic module.

Hence we can define cyclic and periodic cyclic complexes of this cyclic module.
Suppose now that our Hopf algebra $\grch$ is a DG Hopf algebra with the
differential $d$ of degree $1$. Then complexes $\cb$ and
$\cb_{\text{per}}$ have an extra differential defined to be
  $(-1)^nd$ on $\left(\grch\right)^{\otimes n}$
 where we extend $d$ by
\begin{equation}
d(h^1\otimes h^2 \dots \otimes h^n)=
\sum \limits_{i=1}^n (-1)^{\deg h^1 +\dots \deg h^{i-1}}
 h^1\otimes h^2 \dots dh^i \dots\otimes h^n
\end{equation}

We consider the total complexes of the \emph{finite} cochains in the
resulting tricomplexes, and define cyclic and peridic cyclic cohomology
of DG Hopf algebra as cohomology of these complexes. 


 Suppose now that  we are given an action $\pi$ of a differential graded
Hopf algebra $\grch$ on $\W^*$,
which agrees with the differential graded structures on $\grch$ and $\W^*$,
 i.e. in addition to the general properties of Hopf algebra action we have
\begin{align}
\deg \pi(h)(a)&=\deg h+\deg a\\
d\left( \pi(h)(a)\right)&=\pi(dh)(a)+(-1)^{\deg h} \pi(h) (da)
\end{align}
where $h \in \grch$, $a \in \W^*$.
We will often omit $\pi$ from our notations and write just $h(a)$ if it
is clear what action we are talking about.

Suppose that $\dint$ is a closed graded $\s$-trace on $\W^*$,
$\de$-invariant under
the action of $\grch$, i.e.
\begin{align}
\dint \pi(h)(a) b& =\dint a \pi\left(\widetilde{Sh}\right) b \label{defde}\\
\dint ab & =\dint b \pi(\s)(a) \label{defs}
\end{align}
Then one has a  map of cyclic modules
$\chi_{\pi} : (\grch)^{\sharp} \to (\W^*)^{\sharp}$,
given by
\begin{equation} \label{defchi}
\chi_{\pi}(h^1 \otimes h^2 \dots \otimes h^k)(a_0, a_1, \dots a_k)=
\lb \dint a_0\pi(h^1)(a_1) \dots \pi(h^k)(a_k)
\end{equation}
where
\[\lb =(-1)^{  \sum \limits_{ j> i \ge 0} \deg h^j \deg a_i}\]
This map also commutes with the differential $d$, and hence induces a
characteristic map $\chi_{\pi} : \cb(\grch,d) \to \cb(\W^*,d)$, as well
as corresponding maps in cohomology.

\section{Properties of the characteristic map}\label{S:prop}

We will consider now two actions of $\grch$ on $\W^*$ which are
conjugated by the inner automorphism. We will work in the assumption
that $\W^*$ is unital, indicating the changes which need to be made in
the nonunital case in the Remark \ref{nonunit}.
More precisely, let $\rp$ and $\rn$
be two degree-preserving linear maps from $\grch$ to $\W^*$, which
commute with the differentials. We suppose that they are inverse to each
other with respect to convolution:
\begin{equation}\label{dr1}
\sum \rp(\hz)\rn(\ho)=\ve(h)1
\end{equation}
and satisfy cocycle identities:
\begin{align}
\rp(hg)&=\sum \rp(\hz)\pi(\ho)(\rp(g)) \label{dr2}\\
\rn(gh)&=\sum \pi(\hz)(\rn(g)) \rn(\ho) \label{dr3}\\
\rp(1)&=\rn(1)=\rp(\s)=\rn(\s)=1 \label{dr4}
\end{align}
Then one can define a new action $\pi'$ of $\grch$ on $\W^*$ by
\begin{equation}\label{defpi'}
\pi'(h)(a)=\sum (-1)^{\deg \htw \deg a}
\rp(\hz)\pi(\ho)(a)\rn(\htw)
\end{equation}
\begin{lemma}
Equation \eqref{defpi'} defines an action of the DG Hopf algebra $\grch$
on the DG algebra $\W^*$.
\end{lemma}
\begin{proof}
We check all  the required properties.
First
\begin{multline}
\pi'(h)(ab)=\sum(-1)^{\deg ab \deg \htw}
\rp(\hz)\pi(\ho)(ab)\rn(\htw)=\\
\sum(-1)^{(\deg a+\deg b)\deg \hth}(-1)^{\deg a \deg \htw}
 \rp(\hz)\pi(\ho)(a)\pi(\htw)(b)\rn(\hth)=\\
\sum (-1)^{(\deg a+\deg b)\deg \hfi}(-1)^{\deg a (\deg \htw+\deg \hth +\deg
\hf)}\\
\rp(\hz)\pi(\ho)(a)\rn(\htw)\rp{\hth}\pi(\hf)(b)\rn(\hfi)=\\
\sum (-1)^{\deg a \deg \ho}\pi'(\hz)(a) \pi'(\ho)(b)
\end{multline}
Then
\begin{multline}
\pi'(hg)(a)= \\ \sum(-1)^{\deg a (\deg \htw+ \deg \gth)}
\rp(\hz \gz)\pi(\ho \go)(a)\rn(\htw \gtw)=\\\sum (-1)^{\deg a(\deg\hth +\deg
\gtw +\deg \hf)}\\
 \rp(\hz)\pi(\ho)(\rp(\gz)) \pi(\htw \go)(a) \pi(\hth) \rn(\gtw)\rn(\hf)=\\
\pi'(h)\left(\pi'(g)(a)\right) 
\end{multline}

Also
\begin{align}
\pi' (h)(1)=&\sum \rp(\hz) \ho(1) \rn(\hth)=\ve (h)1\\
\pi' (1)(a)=&\rp(1) a \rn(1)=a
\end{align}
and
\begin{multline}
d\left(\pi'(h)(a)\right)=d\sum(-1)^{\deg a \deg \hth}\rp(\hz) \pi(\ho)(a)\rn(\htw)=\\
\sum(-1)^{\deg a \deg \hth}\rp( d\hz) \pi(\ho)(a)\rn(\htw)+\\
\sum(-1)^{\deg a \deg \hth+\deg \hz }\\ \rp(\hz) \left(\pi(d\ho)(a)+
(-1)^{\deg \ho}\pi(\ho)(da)\right)\rn(\htw)+\\
\sum(-1)^{\deg a \deg \hth+\deg\hz +\deg \ho +\deg a }\rp(\hz) \pi(\ho)(a)\rn(d\htw)=\\
\pi'(d h)(a)+
(-1)^{\deg h}\pi'(h)(da)
\end{multline}
\end{proof}

Suppose now that  $\dint$ is the closed $\de$-invariant $\s$-trace for both
actions $\pi$ and $\pi'$. In this case we have two characteristic maps
$\chi_{\pi}$ and $\chi_{\pi'}$ from $\cb(\grch,d)$ to $\cb(\W^*,d)$.
Then we have the following
\begin{prop} \label{inner}
Let  $\pi$ and $\pi'$ be two actions of $\grch$ on $\W^*$, conjugated
by inner automorphisms, and suppose that they both
satisfy conditions \eqref{defde},\eqref{defs}. Let
$\chi_{\pi}$, $\chi_{\pi'}$ be the corresponding characteristic maps.
Then induced maps in cohomology $HC^*\left(\grch,d\right)\to
HC^*\left(\W^* \right)$ coincide.
\end{prop} 
\begin{proof}
Let $M_2(\W^*)=\W^* \otimes M_2(\Cb)$ be the differential graded algebra
of $2\times 2$ matrices over the algebra $\W^*$.We can define an action
$\pi_2$ of $\grch$  on $\W^* \otimes
M_2(\Cb)$ by $\pi_2(h)(a\otimes m)=\pi(h)(a)\otimes m$, where $h\in \grch$, $a\in
\W^*$, $m \in M_2(\Cb)$. Put now
\begin{equation}
\rp_2(h)=\left( \begin{matrix}
                 &\rp(h) &\quad 0\\
                 &0  \quad  &\ve(h)
                 \end{matrix}  
          \right)
          \quad \rn_2(h)=\left( \begin{matrix}
                 &\rn(h) &\quad 0\\
                 &0  \quad  &\ve(h)
                 \end{matrix}  
          \right)          
\end{equation}
It is easy to see that $\rp_2$, $\rn_2$ satisfy equations
\eqref{dr1}-\eqref{dr4}, and hence we can twist the action $\pi_2$ by
$\rp_2$, $\rn_2$ to define a new action $\pi'_2$, as in \eqref{defpi'}.

Consider now the linear functional $\dint_2$ on $M_2(\W^*)$ defined by
\begin{equation}
\dint_2(a\otimes m)= \left(\dint a\right) \left(\tr m\right)
\end{equation} 
Then $\dint_2$ is a closed graded $\de$-invariant $\s$-trace on
$M_2(\W^*)$ with respect to the action $\pi'_2$.
Hence we can define the characteristic map
$\chi_{\pi'_2} : \cb(\grch,d) \to \cb\left(M_2(\W^*), d\right)$

Consider now two imbeddings $i$, $i': \W^* \hookrightarrow M_2(\W^*)$ defined by
\begin{equation}
i(a)=\left(
        \begin{matrix}
        &0 \quad &0\\
        &0 \quad &a
        \end{matrix}
       \right) \quad
i'(a)=\left(
        \begin{matrix}
        &a \quad &0\\
        &0 \quad &0
        \end{matrix}
       \right)
\end{equation}

It is easy to see that $i^* \circ \chi_{\pi'_2} =\chi_{\pi}$ and
$(i')^* \circ \chi_{\pi'_2} =\chi_{\pi'}$.
Now to finish the proof it is enough to recall the well-known fact that $i$ and $i'$
induce the same map in cyclic cohomology. Since we will later need an
explicit
homotopy between $\chi_{\pi}$ and $\chi_{\pi'}$ we give the proof below.

Put $u_t=\exp t \left( \begin{matrix}
                        0& -1&\\
                        1&  0&
                        \end{matrix}
                \right)=
                \left( \begin{matrix}
                        \cos t& -\sin t&\\
                        \sin t&  \cos t&
                        \end{matrix}
                \right)         $
Put $i_t(a)=u_t i(a) u_t^{-1}$. Notice that $i_0=i$, $i_{\pi/2}=i'$.
Consider the family of maps $i_t^*: \cb\left(M_2(\W^*)\right) \to
\cb\left(\W^* \right)$. 
Since we have $\frac{d}{dt}i_t(a) =[g, i_t(a)]$, where $g=\left( \begin{matrix}
                        0& -1&\\
                        1&  0&
                        \end{matrix}
                \right)$
these maps satisfy $\frac{d}{dt}i_t^* =i_t^*L_g$, where $L_g: C^k(M_2(\W^*))
\to C^k(M_2(\W^*))$ is the operator defined by
\[
L_g \phi(x_0,\dots,x_k)=
\sum \limits_{j=0}^k \phi(x_0,\dots,[g,x_j],\dots,x_k)
\]
Define also an operator $I_g :C^k(M_2(\W^*))
\to C^{k-1}(M_2(\W^*))$  by
\begin{equation}\label{defI}
I_g \phi(x_0,\dots,x_{k-1})=
\sum \limits_{j=0}^{k-1} \phi(x_0,\dots,x_j, g,x_{j+1},\dots,x_{k-1})
\end{equation}
Then it is easy to verify that $[b,I_g]=L_g$, $[B,I_g]=0$ and
$[d,I_g]=0$. Hence $L_g=\del I_g+I_g \del$, $\del=\pm d +b +B$.
We conclude that $i_1^*-i_0^*=K\del +\del K$, where the homotopy $K$ is
given by $K\phi=\int_0^{\pi/2} i_t^*I_g$.

Hence
\begin{equation} \label{homotop}
\chi_{\pi'}-\chi_{\pi}=\del H+ H\del
\end{equation}
 where $H=K \circ \chi_{\pi'_2}$
\end{proof}

Now note that the complex
 $\cb\left( \grch \right)$ has a natural weight
filtration by subcomplexes $F^l\cb\left( \grch,d\right)$,
 where
\begin{equation}\label{deffiltr}
F^l\cb\left( \grch,d\right)=\{ \al_1 \otimes \al_2 \dots \otimes \al_j\ |\
\deg \al_1+\deg \al_2 + \dots \deg \al_j \geq l
\end{equation}

Suppose now that $\dint$ has weight $q$, i.e.
\begin{equation}\label{weight}
\dint a=0\text{ if }\deg a \ne q
\end{equation}
Notice that in this case $\chi_{\pi}$ reduces the total degree by $q$
 Then following then proposition is clear:
 \begin{prop} 
 The characteristic map is $0$ on $F^l\cb\left( \grch \right)$ for $l > q$.
 \end{prop}

Let $\cb \left( \grch,d\right)_l$ denote the truncated cyclic bicomplex:
\begin{equation}
\cb \left( \grch,d\right)_l=\cb\left( \grch,d\right)/F^{l+1}\cb\left( \grch,d\right)
\end{equation} Then we immediately have the following
\begin{cor}
The characteristic map $\chi_{\pi}$ defined in
\eqref{defchi} induces the map from the complex
 $\cb \left( \grch,d\right)_q$ to the cyclic complex of the differential
graded algebra $\W^*$.
\end{cor}
This new map will also be denoted $\chi_{\pi}$.
We use the notation 
\begin{equation}\label{deftrun}
HC^*(\grch,d)_l=H^*\left(\cb \left( \grch,d\right)_l\right)
\end{equation}
for the  cohomology of the complex $\cb \left( \grch,d\right)_l$.
The explicit form of the homotopy in the Proposition \ref{inner} now
implies the following
\begin{cor}
Suppose in addition to the conditions of the Proposition \ref{inner}
that \eqref{weight} is satisfied. Then the two maps in cohomology
induced by $\chi_{\pi}, \chi_{\pi'} : \cb \left( \grch,d\right)_q \to
\cb\left(\W^*\right)$ are the same. 
\end{cor}
\begin{proof}
We use the notations of the proof of the Proposition \ref{inner}.
There we established that $\chi_{\pi'}-\chi_{\pi}=\del H+ H\del$. We need only
to verify that $H$ is well defined on the quotient complex $\cb \left( \grch,d\right)_q$.
But since $H=K \circ \chi_{\pi'_2}$, and $\chi_{\pi'_2}$ is
easily seen to be $0$ on  $F^{q+1}\cb\left( \grch,d\right)$, the
result follows.
\end{proof}

\begin{rem} \label{nonunit}
We worked above in the assumption that the DG algebra $\W^*$
is unital. If this is not the case some changes should be made. First of
all, $\rp(h)$, $\rn(h)$ now don't have to be elements of the algebra,
but rather multipliers, such that $\pi'$ defined in \eqref{defpi'} is
a Hopf action. Moreover, we need to require that if $m$ is such a
multiplier, then $\dint ma =\dint am$ $\forall a \in \W^*$. Then the
 the Proposition \ref{inner} remains true. Characteristic maps in this
case take values in the $\cb$ complex of the algebra $\W^*$ with unit
adjoined. The homotopy between two characteristic maps is still given by
explicit formula \eqref{homotop}, which continues to make sense in the
nonunital situation. Indeed, it is sufficient to define $I_g
\chi_{\pi'_2}$, which can be defined by the same formula as above,
provided one  treats $g$ and $\pi'_2(h)(g)$ as
multipliers of the algebra $M_2\left(\W^*\right)$, with
$\pi'_2(g)$ defined to be
$ \left(
\begin{matrix}
&0 &\rp(h) \\
&-\rn(h) &0
\end{matrix}
\right)$
\end{rem}
Finally, we collect all the information we will need to use in the next sections.
\begin{thm}
Let $\left(\W^*,d \right)$ be a differential graded algebra, and
$\dint$ a linear functional on $\W^*$ of weight $q$, and let $\ca=\W^0$
be the degree $0$ part of $\W^*$.
 Let $\pi$ be an action of the DG Hopf algebra $\grch$ act on the DGA $\W^*$.
 Suppose that $\dint$ is a $\de$-invariant $\s$-trace with respect to $\pi$.
 Then
characteristic map \eqref{defchi} defines a map in cohomology $HC^i(\grch)_q \to
HP^{i-q}(\ca)$. Suppose now that $\pi'$ is another action of $\grch$ on
$\W^*$, obtained from $\pi$ by twisting by a cocycle \eqref{defpi'}.
Then if $\dint$ is a   $\de$-invariant $\s$-trace with respect to
$\pi'$, the maps in cohomology $HC^i(\grch)_q \to
HP^{i-q}(\ca)$ induced by $\chi_{\pi}$, $\chi_{\pi'}$ are the same.
\end{thm}


\section{Secondary characteristic classes}\label{S:forms}
Let $M$ be a manifold, and let
$\G$ be a discrete pseudogroup of diffeomorphisms of $M$, acting from
the right.

By this we mean a set $\G$ such that every element $g$ of $\G$
defines a  local diffeomorphism 
of $M$, i.e. diffeomorphism $g: \Dom g \to \Ran g$, where
$\Dom g$, $\Ran g \subset M$ -- open subsets of $M$, and  that
we have partially defined operations of composition and inverse such
that 
\begin{enumerate}
\item If $g \in \G$ then $g^{-1} : \Ran g \to \Dom g$ is also in $\G$.
\item If $g_1, g_2 \in \G$ then $g_1g_2$ with domain
$g_1^{-1}\left(\Dom g_2\cap \Ran g_1 \right)$ and range $g_2
\left(\Dom g_2\cap \Ran g_1 \right)$ is in $\G$.
\item $id :M \to M$ is in $\G$.
\end{enumerate}
Note that we use a wide definition of pseudogroups, and do not include
any saturation axioms.

Let  $E$ be a trivial vector bundle on $M$, equivariant with respect
to the action of $\G$. In other words,  every $g \in \G$ defines for
every $x \in \Dom g$ a linear map $E_x \to E_{xg}$. 

For the rest of the paper we suppose the following:

\emph{If $g_1$ and $g_2 \in \G$ are such that they induce the same
diffeomorphisms and the same action on the bundle, then $g_1= g_2$}.

With this data one can associate the following groupoid $\cg$:
the objects are the points of $M$ and the morphisms $x \to y$, $x$, $y \in M$
are given by $g \in \G$ such that $g(x)=y$, with the composition given
by the product in $\G$.
Let $\ca$ denote the convolution algebra of this groupoid, i.e. the
cross-product $C_0^{\ify}(M) \rtimes \G$. Let
$\W^* =\left(\W^*(M)\rtimes \G, d\right)$ denote the differential
graded algebra of forms on $\cg$ with the convolution product, where
the differential $d$ is the de Rham differential.
We will use the usual cross-product notations $\w U_g$ for the elements
of this algebra, where $\w\in \W^*(M)$, $g \in \G$. Since $\G$ is, in
general  a pseudogroup, we suppose that 
\begin{equation}
\supp \w \subset \Dom g
\end{equation}

Fix a trivialization of $E$. The action of $\G$ on the bundle defines
then a homomorphism
\begin{equation}\label{defh}
h: \cg \to GL_n(\Rb)
\end{equation}

Let $\chg$ denote the differential graded Hopf algebra of the forms on
$GL_n(\Rb)$, with the product given exterior multiplication,
coproduct, antipode and counit induced respectively
by the product $ GL_n(\Rb) \times 
GL_n(\Rb) \to  GL_n(\Rb)$, inverse $ GL_n(\Rb) \to  GL_n(\Rb)$ and the
inclusion $1 \to  GL_n(\Rb)$. The differential is given by the de Rham
differential on forms.

We now show that the map \eqref{defh} allows one to define an action of
$\chg$ on $\W^*$. 
\begin{prop}
The map $\chg \otimes \W^* \to \W^*$ given by
\begin{equation}
\pi(\al) (\w) =h^*(\al)  \w
\end{equation}
where $\al \in \chg$, $\w \in \W^*$
defines an action of the differential graded Hopf algebra
$\chg$ on the differential graded algebra $\W^*$.
\end{prop}
\begin{proof}
We have:
\begin{equation}
\pi(\al_1 \al_2) (\w)=h^*(\al_1 \al_2)  \w =
h^*(\al_1)  h^*( \al_2)  \w=\pi(\al_1) \left(\pi(\al_2) (\w)\right)
\end{equation}

Next, if we write $\De \al = \sum \limits_k \al_{(0)}\otimes \al_{(1)}$
we have:

\begin{multline}
\pi(\al)(\w_0 \w_1)(g)=h^*(\al)(g)  \w_0   \w_1(g)=
h^*(\al)(g) \sum \limits_{g_0g_1=g}\w_0(g_0) \w_1^{g_0}(g_1)=\\
\sum \limits_{g_0g_1=g}
\sum \limits_k h^*(\al_{(0)})(g_0)
h^*(\al_{(1)})^{g_0}(g_1) \w_0(g_0) \w_1^{g_0}(g_1)=\\
\sum \limits_{g_0g_1=g}
\sum \limits_k (-1)^{\deg \w_0 \deg \al_{(0)}}
h^*(\al_{(0)})(g_0) 
\w_0(g_0)   h^*(\al_{(2)})(g_1) \w_1^{g_0}(g_1)=\\
\sum \limits_k (-1)^{\deg \w_0 \deg \al_{(0)}}
\pi(\al_{(0)})(\w_0) \pi(\al_{(1)})(\w_1)
\end{multline}
Also, if $M$ is compact the algebra $\W^*$ has a unit  given by the function
\begin{equation}
e(g)=
\begin{cases}
1 &\text{ if $g$ is a unit}\\
0 &\text{ otherwise}
\end{cases}
\end{equation}
Then
\begin{equation}
\pi(\al)(e)=h^*(\al) e =\e(\al) e
\end{equation}
Finally, we have
\begin{multline}
d(\pi(\al(\w))) =
d\left(h^*(\al) \w)\right)=\\
 h^*(d \al) \w +(-1)^{\deg \al} h^*(\al) d\w
 =\pi((d\al))(\w) +(-1)^{\deg \al} \pi(\al) (d\w)
\end{multline}
\end{proof}

We now have a natural inclusion $i: M \hookrightarrow \cg$ as the
space of units. Then we define a graded trace $\dint$ on $\W^*$ by
\begin{equation}
\dint \w =\int \limits_M i^* \w
\end{equation}
\begin{prop}
The graded trace $\dint$ is closed under the de Rham differential and is
invariant under the action of $\ch$, i.e.
\begin{align}
&\dint d\w=0\\
&\dint \al(\w) =\e(\al) \dint \w
\end{align}
\end{prop}
\begin{proof}
The first identity is clear, the second follows from the fact that
$h \circ i : M \to GL_n(\Rb)$ is a constant map, taking the value 1.
\end{proof}

Hence we have a map $\cb_q\left(\chg,d\right) \to  \cb\left(\W^*,d\right)$
where $q=\dim M$, which also gives us a map
\begin{equation}\label{chi}
\chi:HC_q^*\left(\chg,d\right) \to HP^*\left (C_0^{\ify}(M)\rtimes \G\right)
\end{equation} 
Definition of the action of $\chg$ on $\W^*$, and hence definitions
of the map given by
\eqref{chi} apriori
depends on the choice of
trivialization of $E$,
but we will now show that this
map  is independent of the choice of trivialization.

\begin{prop}
Suppose we use another trivialization of $E$ to define an action of
$\chg$ on $\W^*$. Then the two actions are conjugated by the inner
automorphisms. 
\end{prop}
\begin{proof}
Let us chose another trivialization of the bundle $E$, and let $U(x)$ $x
\in M$
be a transition matrix between the two bases of the fiber $E_x$.
Then we have a new map $h': \cg \to GL_n(\Rb)$, related to $h$ by
\begin{equation}
h'(\g)= U(s(\g))h(\g)U^{-1}(r(\g))
\end{equation}
Let $\pi'$ denote the action corresponding to the map $h'$.

Consider now the pull-back $U^*: \W^*(GL_n(\Rb))\to \W^*(M)$ as a map
$\chg \to \W^*$, where we consider forms on $M$ as the form on $\cg$
which is $0$ outside the space of units. When $M$ is not compact, we
obtain not an element in algebra, but rather a multiplier. Hence we see
that if  we   define
\begin{equation}
\rp(\al) =U^*(\al)
\end{equation}
and
\begin{equation}
\rn(\al)=(U^{-1})^*(\al) =\rp(S\al)
\end{equation}
we will have
\begin{equation}
\pi'(\al)(\w)=\sum (-1)^{\deg \al_{(2)} \deg \w}
\rp(\al_{(0)})\pi(\al_{(1)})(\w)\rn(\al_{(2)})
\end{equation}
\end{proof}

We can now summarize the results as follows.
\begin{thm}
Let $\G$ be a discrete pseudogroup acting on the manifold $M$
of dimension $q$ by
orientation preserving diffeomorphisms. Let $E$ be a $\G$-equivariant
trivial bundle of rank $n$ on $M$. Let $\chg$ be the DG Hopf algebra of the
differential forms on the Lie group $GL_n(\Rb)$. Then we have a map
\begin{equation}
\chi: HC^i_q(\chg, d) \to HP^{i-q}\left(C_0^{\ify}(M)\rtimes \G\right)
\end{equation}
which is independent of the trivialization of $E$. In conjunction with
the equation \eqref{cohom} below it gives a map
\begin{equation}
\bigoplus \limits_{m \in \Zb} H^{i-2m}\left (W(\mathfrak{gl_n},
O_n)_q\right) \to HP^{i-q}\left(C_0^{\ify}(M)\rtimes \G\right)
\end{equation}
\end{thm}

\section{Relation with Weil algebras}\label{S:rel}
In this section we use methods of \cite{kt75}, \cite{shs78} and \cite{dup76,dup78}
to identify the cyclic
cohomology $HC^*(\chg, d)_q$. It turns out that  
\begin{equation}\label{cohom}
HC^i(\chg, d)_q =\bigoplus \limits_{m\ge0} H^{i-2m}\left (W(\mathfrak{gl_n},
O_n)_q\right).
\end{equation}
where $H^*\left (W(\mathfrak{gl_n},
O_n)_q\right)$  is  the cohomology of truncated Weil algebra (cf. \cite{kt75}).
As a matter of fact, one can work with the DG
Hopf algebra $\ch(G)$ of differential forms on any almost connected Lie
group, and the result in this case is
\begin{equation}
HC^i(\ch(G), d)_q =\bigoplus \limits_{m\ge0} H^{i-2m}\left (W(\mathfrak{g},
K)_q\right).
\end{equation}
where $K$ is the maximal compact subgroup of $G$.
 Computation of the Hochschild cohomology of this Hopf algebra
is essentially contained in \cite{kt75}, and with little care using ideas
from \cite{dup76,dup78} one recovers the cyclic cohomology.

Let $G$ be a Lie group with finitely many
connected components. Similarly to the previous section we can
define a differential graded Hopf algebra $\ch(G)$. Let $K$ be the
maximal compact subgroup of $G$. We will now construct the map of complexes
from the truncated relative Weil algebra $W(\mathfrak{g}, K)_q$ to the
complex $\cb_q \left(\ch(G)\right)$.

Let $NG$ denote the simplicial manifold  with
$NG_p=\underbrace{G\times G\times \dots  G}_p$
The simplicial structure  is given by 
the face maps  
\begin{equation}
\del_i(g_1, g_2,\dots,g_k) =
\begin{cases}
(g_2,\dots,g_k) &\text{ if } i=0\\
(g_1, g_2,\dots
g_ig_{i+1},\dots,g_k) &\text{ if } 1\le i \le k-1\\
(g_1, \dots,g_{k-1}) &\text{ if } i= k
\end{cases}
\end{equation}
 and degeneracy maps 
 \begin{equation}
\s_i(g_1, g_2,\dots,g_k) =
(g_1,  \dots, g_{i-1}, 1, g_i,\dots, g_k)
\end{equation}
The geometric realization of this simplicial manifold is the classifying
space $\bg$. It is a union of manifolds $\De^p \times NG_p$ with the
modulo the equivalence relation (cf. \cite{dup78}).

We will also consider  simplicial manifold $\bar{N}G$, with
$\bar{N}G_p =\underbrace{G\times G \times \dots \times G}_{p+1}$.
The face and degeneracy maps are given by
\begin{equation}
\del_i(g_0,g_1, g_2,\dots,g_k)=(g_0,\ldots \hat{g}_i, \ldots, g_k)
\end{equation}
\begin{equation}
\s_i(g_0,g_1, g_2,\dots,g_k) =(g_0,  \ldots, g_{i-1},g_i, g_i, g_{i+1}
,g_k)
\end{equation}
The geometric realization of this simplicial manifold
is $\eg$.
The map $pr : \bar{N}G \to NG$ given by
\begin{equation}
pr (g_0, g_1, \dots, g_p)=(g_0g_1^{-1}, g_1g_2^{-1}, \dots, g_{p-1}g_p^{-1})
\end{equation}
defines a simplicial principal $G$-bundle $\eg \to \bg$.
Simplicial manifolds $\bar{N}G$ and $NG$ moreover have a cyclic
structure, i.e. an action of  the cyclic groups $\Zb_{p+1}$ 
on the $p$-th component, which satisfy all the necessary relations with
the face and degeneracy maps.
The actions are given on $\bar{N}G$ by
\begin{equation}
\tau_p(g_0, g_1, \dots,g_p)=(g_1, g_2, \dots, g_p,g_0)
\end{equation}
Since the maps $\tau_p$ are $G$-equivariant, they induce corresponding
actions on $NG$:
\begin{equation}
\tau_p(g_1, g_2,\dots,g_p)=(g_2, g_3, \dots, g_p,(g_1g_2\dots g_p)^{-1})
\end{equation}
We will identify the $p$-cochains for the Hopf algebra $\ch(G)$ with the
forms on $NG_p$. Under this identification the
simplicial structure on the Hopf cochains corresponds to the one
induced by the simplicial structure on $NG$, the de Rham differential on
the Hopf cochains corresponds to the de Rham
differential on $NG$, and cyclic structure on the Hopf cochains is
induced by the cyclic structure on $NG$.
Filtration by the form degree on the Hopf algebra cochains corresponds
to the filtration by the form degree on the 
manifold $NG$.

We will now construct the map $\mu$ from the complex $W(\mathfrak{g},K)$
to the simplicial-de Rham complex of forms on $NG$, which preserves
filtration on these complexes.
We do it by constructing the map from $W(\mathfrak{g},K)$ to the
complex of simplicial forms on $\bg$, and then applying the integration map.
The complex of simplicial forms on $\bg$ has a natural bigrading. 
Let $\te$ be the Maurer-Cartan form on $G$. Let $p_i :\bar{N}G=G^{p+1} \to
G$ be the projection on $i$-th component. Consider on $\eg_p$ the
$\mathfrak{g}$ valued differential form $\w \sum t_i  \te_i$, where
$\te_i=p_i^*\te$. It
defines a simplicial connection in the bundle $\eg \to \bg$.
The standard construction defines a differential graded algebra
homomorphism $\psi$ from $W(\mathfrak{g},K)$ to the complex of $K$-basic simplicial forms on
$\eg$, which we identify with forms on the space $\eg/K$.
This space is a bundle over $\bg$ with the fiber $G/K$.
This bundle has a section, which can be explicitly described as follows.
Since $G/K$ has a natural structure of a manifold of constant negative
curvature, for any finite set of points $x_0$, $x_1$, ..., $x_k \in G/K$
one can construct a canonical simplex in $G/K$, i.e. a map
$\s(x_0, x_1, \dots, x_k):\De^k \to G/K$,
 with vertices $x_0$, $x_1$, ..., $x_k$, and this
construction agrees with taking faces of a simplex, and is
$G$-equivariant:
\begin{equation}
\s(gx_0, gx_1, \dots, gx_k)(t_0,t_1, \dots, t_k)=g\s(x_0, x_1, \dots, x_k)(t_0,t_1, \dots, t_k)
\end{equation}
Denote by $\pi$ the canonical projection $G \to G/K$.
 Then the section $s$ is given by the simplicial map defined by the
following formula, where we write just $\s$ for $\s(\pi(1),\pi(g_1), \pi(g_1g_2), \dots
\pi(g_1\dots g_k)(t_0,t_1, \dots,t_k)$: 
\begin{multline} s(g_1,g_2, \dots, g_k;t_0,t_1, \dots t_k)= \\
\bigl( \s^{-1},
 g_1\s^{-1},
 g_1g_2\s^{-1},  \dots
g_1\dots g_k\s^{-1}; t_0, t_1, \dots t_k\bigr)
\end{multline}
Notice that this section intertwines   the actions of the cyclic group
on the spaces $\bg$ and $\eg/K$.

Consider now the map $s^* \psi$. It is clearly a homomorphism from the
differential graded algebra $W(\mathfrak{g},K)$ to the differential graded
algebra of simplicial forms on $\bg$. We will show now that it preserves
filtration's on both algebras.
First we need the following statement:
\begin{lemma} \label{type}
If $\xi$ is a horizontal form on $\eg$ of the type $(k,l)$,then
 $s^*\xi$ is also of the type $(k,l)$.
\end{lemma}
\begin{proof}
Since $\xi$ is horizontal, it can be written as a sum of the expressions
of the form
$f pr^* \zeta$, where $f$ is a function,  $pr:\eg \to \bg$ --
projection, and $\zeta$ is a form on $\bg$ of the type $(k,l)$.
Then $s^*\left( f pr^* \zeta \right)=(s^*f) \zeta$ is also of the type   $(k,l)$.
\end{proof}
\begin{prop}
 The homomorphism $s^* \psi$ agrees with
filtrations on the Weil algebra and on the forms on $\bg$.  
\end{prop}
\begin{proof}
The curvature of the connection $\w$ is a horizontal form $\W$ on $\bg$,
given by
\begin{equation}
\W =\sum dt_i \te_i +\sum t_i d\te_i +\sum \limits_{i<j}t_it_j[\te_i,\te_j]
\end{equation}
Hence $\W$ has only components of the type $(1,1)$ and $(0,2)$. The
statement of the lemma will then follow from the fact that $s^*\W$ 
also has only components of the type $(1,1)$ and $(0,2)$.
But this follows from  the  Lemma \ref{type}.
\end{proof}

We can now apply the integration map and obtain the map $\mu$ from the
Weil algebra to the simplicial-de Rham complex of $NG$. Since the
integration map respects  filtrations, the resulting map $\mu$ also
respects filtrations. We identify the Hochschild complex of $\ch(G)$
with the simplicial-de Rham complex of $NG$. Results of \cite{kt75},
\cite{shs78} imply that this is actually an isomorphism, i.e. 
\begin{equation}
HH^i(\ch(G), d)_q = H^{i}\left (W(\mathfrak{g},
K)_q\right)
\end{equation}

But since the connection
$\w$ and the section $s$ are invariant under the cyclic action, the
resulting Hochschild cochains are actually cyclic. 
This implies that  Connes' long exact sequence is equivalent to the
collection of short exact sequnces 
\begin{multline}
0\to HC^{i-2}(\ch(G), d)_q \overset{S}\to HC^i(\ch(G), d)_q \overset{I}
 \to \\
 \overset{I} \to HH^i(\ch(G), d)_q \to  0
\end{multline}
and the map $I$ splits. Hence
\begin{multline}
HC^i(\ch(G), d)_q = \\ \bigoplus \limits_{m\ge0}HH^{i-2m}(\ch(G), d)_q=
\bigoplus \limits_{m\ge0} H^{i-2m}\left (W(\mathfrak{g},
K)_q\right).
\end{multline}

Explicitly, maps $H^{i-2m}\left (W(\mathfrak{g},
K)_q\right) \to HC^i(\ch(G), d)_q$ are given 
by $S^m \circ \mu$, where we consider $\mu$ as a map 
into cyclic complex.

\section{Relation with other constructions}\label{S:other}
Suppose, as before, that we have an orientation-preserving action of
discrete group $\G$ on an oriented manifold $M$, and an equivariant
trivial bundle $E$ over $M$.
Then results of previous sections provide us a map $H^*(W(\mathfrak{g},O_n))
\to HP^*(\ca)$, where $\ca=C_0^{\ify}(M)\rtimes \G$.
We also have a construction of the map $H^*(W(\mathfrak{g},O_n)) \to
H^*(M_{\G})$ (see e.g. \cite{kt75,br,bh})
 where $M_{\G}=M\times_{\G}\rm{E}\G$ is the homotopy quotient. In this
section we prove that these constructions are compatible, i.e. that
the following diagram is commutative
\begin{equation}\label{coin}
\xymatrix{&H^*(W(\mathfrak{g},O_n)) \ar[rr] \ar[rrd] & &H^*(M_{\G})\\
          &                                         & &HP^*(\ca)\ar[u]_{\Phi}
         }
\end{equation}
where $\Phi$ is the canonical map given by Connes \cite{cn86, cn94}.

The proof goes as follows. We construct a map $\Psi$ from the complex
computing $H^*(M_{\G})$ to the cyclic complex $\cb(\W^*,d)$, where
$\W^*=\W^*_0(M)\rtimes \G$, which has the following properties.
First, it agrees with the map $\Phi$, in the sense that the following
diagram is commutative:
\begin{equation}
\xymatrix{&H^*(M_{\G}) \ar[rr]^{\Phi} \ar[rrd]^{\Psi} & &HP^*(\ca)\\
          &                                         & &HP^*(\W^*,d)\ar[u]_{R}
         }
\end{equation}  
where the map $R$ is defined by \eqref{defR}.
Then it is clear from the definitions that
 the diagram similar to \eqref{coin} is valid with the map $\Psi$
 already on the level of cochains, not just cohomology.

 The definition of the map $\Psi$ is the following.
Recall that the cohomology of $M_{\G}$ can be computed by the following
bicomplex $\cc^{*,*}$. $\cc^{k,l}$ denotes the set of totally antisymmetric
functions $\tau$ on $\underbrace{\G\times\G \dots \times \G}_{k+1}$ with values in
$-l$-currents on $\Dom g_0 \cap \Dom g_1 \dots \cap \Dom  g_k$, which
satisfy the invariance condition 
\begin{equation}
\tau(gg_0,gg_1,\dots,gg_{k})=\tau(g_0,g_1,\dots,g_k)^{g^{-1}}.
\end{equation}

The two differentials of this complex are given by the 
the group cohomology complex differential given on $\cc^{k,l}$ by
\begin{equation}
(d_1 \tau) (g_0,g_1,\dots,g_{k}, g_{k+1})= (-1)^l\sum \limits_{j=0}^{k+1}
(-1)^j \tau (g_0,g_1, \dots, \hat{g}_j, \dots, g_{k+1})
\end{equation}
and the de Rham
differential $d$ given by
\begin{equation}
(d_2\tau) (g_0,g_1\dots,g_k)=d\left(\tau(g_0,g_1, \dots,g_k)\right)
\end{equation}

We now define the map $\Psi$ from the complex
 $\cc^{*,*}$ to the cyclic complex $\cb(\W^*,d)$, where
$\W^*=\W^*_0(M)\rtimes \G$, by the following formula.

\begin{multline}
\Psi (\tau)( \w_0U_{g_0}, \w_1U_{g_1},\dots, \w_kU_{g_k})=\\
\begin{cases}
(-1)^{kl}
\langle \tau(1,g_0,g_0g_1,\dots,g_0\dots g_{k-1}),
\w_0\w_1^{g_0}\dots \w_k^{g_0\dots g_{k-1}}\rangle &\text{ if
  }g_0\dots g_k=1\\
0 &\text{ otherwise}
\end{cases}
\end{multline}

Map $\Psi$ satisfies the following identities
\begin{align}
b\Psi(\tau)&=\Psi(d_1\tau)\\
d\Psi(\tau)&=\Psi(d_2\tau)\\
B\Psi(\tau)&=0
\end{align}
and hence it is indeed a map of complexes.
It is clear from the definition of the map $\Psi$ that the diagram
obtained from the diagram \eqref{coin} by replacing $\Phi$ by $\Psi$
commutes, even on the level of complexes. 
It remains to prove that the map $\Psi$ induces the same map in cohomology as the
map $\Phi$.

To do this we note that
Connes' map $\Phi$ is characterized uniquely by the following 
properties (cf. \cite{gor99a,gor99b}):

\begin{enumerate}
 \item 
 $\Phi$ takes values in the part of cyclic cohomology supported at
 the identity of the group. 
 \item 
 Let $M$ be an  oriented manifold on which discrete group $\G$ acts freely and
 properly preserving orientation. Then any class in $c\in H^*(M_{\G})$ can be represented by
   a $\G$-invariant current $C$ on $M$. Then 
$\Phi(c)$ is the same as the class of the cyclic cocycle defined by
\begin{multline}
\phi(a_0U_{g_0}, a_1U_{g_1}, \dots, a_kU_{g_k})=\\
\begin{cases}
\langle C , a_0 da_1^{g_0} \dots da_k^{g_0g_1\dots g_{k-1}}\rangle &\text{ if
  } g_0g_1\dots g_k=1 \\
0 &\text{ otherwise}
\end{cases}  
\end{multline}

\item
Let $X$ be another oriented $\G$-manifold. Then  one has the following
commutative diagram
\begin{equation}
\xymatrix{
&H^*\left( (M\times X)_{\G}\right)\ar[r]^-{\Phi} &HP^*\left(C_0^{\ify}(M\times
X)\rtimes \G\right)\\
&H^*\left( M_{\G}\right)\ar[u] \ar[r]^-{\Phi} &HP^*\left(C_0^{\ify}(M)\rtimes \G \right)\ar[u]
}
\end{equation}
The left vertical arrow here is induced by the natural map $(M\times
X)_{\G} \to M_{\G}$ and the right one is induced by the product with
the transverse fundamental class of $X$.

It is easy to see that the map $\Psi$ satisfies the same properties,
and hence gives the same map in cohomology.
\end{enumerate}


\end{document}